\documentclass[11pt]{article}
\usepackage{amsfonts,latexsym,amssymb,epsfig}
\usepackage{amsmath,amsthm,amstext,amscd}
\usepackage{fullpage}
\usepackage{euscript}
\begin{document}
\bibliographystyle{plain}
\newtheorem{theorem}{Theorem}
\newtheorem{proposition}[theorem]{Proposition}
\newtheorem{corollary}[theorem]{Corollary}
\newtheorem{note}[theorem]{Note}
\newtheorem{lemma}[theorem]{Lemma}
\newtheorem{definition}[theorem]{Definition}
\newtheorem{observation}[theorem]{Observation}
\newcounter{fignum}
\newcommand{\figlabel}[1]
           {\\Figure \refstepcounter{fignum}\arabic{fignum}\label{#1}}
\newcommand{\ignore}[1]{}
\def\F2{{\{0,1\}}}
\def\eps{{\epsilon}}
\def\tP{\tilde{P}}
\def\ttP{\tilde{\tilde{P}}}
\def\rhoi{\rho^{-1}}
\def\hM{\hat{M}}
\def\hpM{\hat{M'}}
\def\halpha{\hat{\alpha}}
\def\hpalpha{\hat{\alpha}'}
\def\tO{\tilde{O}}
\def\tOmega{\tilde{\Omega}}
\def\S{{\Sigma}}
\def\hn{{\lfloor n/2\rfloor}}
\def\ox{{\overline{x}}}
\def\I{\EuScript{I}}
\def\cardI{{I}}
\newcommand{\rank} {\mbox {rank}}
\newcommand{\schreier} {\mbox {sc}}
\newcommand{\gap} {\mbox {gap}}
\title{Simple Permutations Mix Even Better}
\author{
\parbox{8cm}{\centering
Alex Brodsky\\
Department of Computer Science\\
University of Toronto\\
abrodsky@cs.toronto.edu
}
\parbox{8cm}{\centering 
Shlomo Hoory\footnote{
Research is supported in part by an NSERC grant and a PIMS postdoctoral 
fellowship.}\\
Department of Computer Science\\
University of British Columbia\\
shlomoh@cs.ubc.ca}
}
\maketitle

\begin{abstract}
We study the random composition of a small family of $O(n^3)$ 
simple permutations on $\{0,1\}^n$. 
Specifically we ask how many randomly selected simple permutations need be 
composed to yield a permutation that is close to $k$-wise independent.
We improve on the results of Gowers~\cite{Go96} and 
Hoory et al.~\cite{HMMR04} and show that up to a polylogarithmic factor, 
$n^2 k^2$ compositions of random permutations from this family suffice.
In addition, our results give an explicit construction of a degree 
$O(n^3)$ Cayley graph of the alternating group of $2^n$ objects with 
a spectral gap $\Omega(2^{-n}/n^2)$, which is a substantial improvement over
previous constructions.

\end{abstract}



A naturally occurring question in cryptography is how well the composition 
of simple permutations drawn from a simple distribution resembles a random 
permutation.
Although such constructions are a common source of security for 
block ciphers like DES and AES, 
their mathematical justification (or lack thereof) is troubling.

This  motivated the investigation of Hoory et al.~\cite{HMMR04} who considered
the notion of almost 
{\em $k$-wise independence}. Namely, that the distribution obtained when 
applying a permutation from a given distribution to any $k$ distinct 
elements is almost indistinguishable from the distribution obtained when 
applying a truly random permutation.
Therefore, the question is
how close is the composition of $T$ random simple permutations  
to $k$-wise independent?  

Another motivation is a fundamental open problem in the theory of expanding
graphs. 
\footnote{A solution to this problem was announced recently by 
Kassabov~\cite{Ka05}.}
Namely, the problem of constructing a constant degree expanding 
Cayley graph of the symmetric group. 
A possible relaxation of this problem is to  ask whether one can find a small 
set of simple permutations such that its action on $k$ points yields an 
expanding graph. 

It turns out that these two problems reduce to bounding the mixing time and 
the spectral gap of the random walk on the {\em same} graph.
This walk, $P$, is defined on the state space 
of $k$-tuples of distinct elements from the $n$-dimensional binary cube. 
In each step it randomly selects a simple permutation and
applies it to each of the $k$ elements at its current position.
The mixing time, $\tau(\eps)$, is the number of steps needed to 
come $\eps$-close to the uniform distribution (in total variation distance),
and the spectral gap, $\gap(P)$, is the difference between
the two largest eigenvalues of $P$'s transition matrix.

Following the construction of DES, and previous work by Gowers~\cite{Go96} 
and Hoory et al~\cite{HMMR04}, 
we consider the class of {\em width $2$ simple permutation}, denoted $\S$.
The action of such a permutation on an element of the $n$-dimensional 
binary cube is to XOR a single coordinate with a Boolean function of $2$
other coordinates; there are $16n(n-1)(n-2)$ such 
permutations.  

These problems were first considered by Gowers~\cite{Go96} who gave an
$\tO(n^3 k (n^2+k)(n^3+k) )$
\footnote{Notation $\tO$ suppresses a polylogarithmic factor
in $n$ and $k$.} 
 bound on the mixing time, by lower bounding the 
spectral gap $1/\gap(P) = \tO(n^2 (n^2+k)(n^3+k) )$.
Subsequently, Hoory et al.~\cite{HMMR04} improved the bound on the mixing 
time to $\tO(n^3 k^3 )$ by proving that $1/\gap(P) = \tO(n^2 k^2)$.
Both results were achieved using the {\em canonical paths} technique,
and neither result applies for $k > 2^{n/2}$.
Using the comparison technique, in conjunction with the theory of reversible
computation, we give better bounds for all values of $k$ 
up to the largest conceivable value, $k=2^n-2$.

\begin{theorem}\label{k2n2:theorem}
$\tau(\eps) = \tO( n^2 k^2 \cdot \log(1/\eps) )$,
as long as $k \leq 2^{n/50}$.
\end{theorem}

\begin{theorem}\label{k2n3:theorem}
$1/\gap(P) = O( n^2 k )$
for all $k \leq 2^n-2$.
\end{theorem}

Using the well known connection between the mixing time and the
spectral gap Theorem~\ref{k2n3:theorem} implies:

\begin{corollary}
$\tau(\eps) = O( n^2 k \cdot (nk+\log(1/\eps)) )$
for all $k \leq 2^n-2$.
\end{corollary}

The proofs of both Theorems are based on the comparison technique for 
Markov chains~\cite{DiSa93}.  
To prove Theorem~\ref{k2n3:theorem} we compare the random walk $P$
either to a Glauber dynamics Markov chain or to the random walk on
the alternating group using $3$-cycles. 
To prove Theorem~\ref{k2n2:theorem} we observe that 
after a short preamble the random walk $P$ is almost surely in 
a ``generic'' state.
Consequently, it suffices to bound the mixing time of a Markov chain
restricted to ``generic'' states.  
To this end we again employ the comparison technique, 
but with a better comparison constant.
In all cases we construct the multicommodity flows required 
by the comparison technique using ideas from the 
theory of reversible computation.


It follows from \cite{HMMR04,MaPi04} that these results
apply also in the more general setting of adaptive adversaries 
(see references for a definition).
\begin{corollary}\label{strongk2n2:corollary}
Let $T$ be the minimal number of random 
compositions of independent and uniformly distributed permutations from 
$\S$ needed to generate a permutation which is 
$\eps$-close to $k$-wise independent against an adaptive adversary.
Then $T = \tO( n^2 k^2 \cdot \log(1/\eps) )$ for $k \leq 2^{n/50}$,
and $T = O( n^2 k \cdot (nk + \log(1/\eps)) )$ for $k \leq 2^n-2$.
\end{corollary}


\section{Preliminaries}

Let $f$ be a random permutation on some base set $X$. 
Denote by $X^{(k)}$ the set of all $k$-tuples of distinct elements from $X$.
We say that $f$ is $\eps$-close to $k$-wise independent if for every
$(x_1,\ldots,x_k) \in X^{(k)}$ the distribution of 
$(f(x_1),\ldots,f(x_k))$ 
is $\eps$-close to the uniform distribution on $X^{(k)}$.
We measure the distance between two
probability distributions $p, q$ by the total variation distance, defined by
\begin{eqnarray*}
d(p,q) = \frac 1 2 ||p-q||_1 = \frac 1 2 \sum_{\omega} |p(\omega)-q(\omega)|
= \max_A \sum_{\omega \in A} p(\omega)-q(\omega).
\end{eqnarray*}

Assume a group $H$ is acting on a set $X$ and let $S$ be a subset of $H$
closed under inversion. Then the {\em Schreier graph} 
$G=\schreier(S,X)$ is defined 
by $V(G)=X$ and $E(G)= \{(x,xs): x \in X,\, s \in S\}$.
For a sequence $\omega=(s_1,\ldots,s_\ell) \in S^l$ we denote
$x \omega = x s_1 \cdots s_\ell$, and we sometimes refer by $x \omega$ to
the walk $x, x s_1, \ldots, x s_1 \cdots s_\ell$.

The {\em random walk} $X_0,X_1,\ldots$
associated with a $d$-regular graph $G$ is defined by the 
transition matrix $P_{vu} = \Pr[X_{i+1}=u|X_i=v]$ which is $1/d$ if 
$(v,u) \in E(G)$ and zero otherwise. 
The uniform distribution $\pi$ is stationary for this Markov chain.
If $G$ is connected and not bipartite, we know that given any initial 
distribution of $X_0$, the distribution of $X_t$ tends to the uniform 
distribution. The mixing time of $G$ is
$\tau(\eps) = \max_{v \in V(G)} \min \{ t : d(P^{(t)}(v,\cdot),\pi) < \eps \}$,
where $P^{(t)}(v,.)$ is the probability distribution of $X_t$ given that
$X_0=v$.
It is not hard to prove (see~\cite[Lemma 20]{AlFi}) that
\begin{eqnarray}\label{submultmix}
\tau(2^{-\ell-1}) \leq \ell \cdot \tau(1/4).
\end{eqnarray}
Let $1=\beta_0 \geq \beta_1 \geq \cdots \geq \beta_{|V(G)|}$ be the eigenvalues
of the transition matrix $P$.
We say that this random walk is lazy if for some constant $\delta>0$ we have
$P_{vv} \geq \delta$ for all $v \in V(G)$. 
We denote the spectral gap $1-\beta_1$ of the Markov chain $P$ by $\gap(P)$.

Two fundamental results relating the spectral gap of a Markov chain to
its mixing time are the following:
\begin{theorem}\label{gapmix:theorem}(\cite[Proposition 3]{DiSt91})
If the random walk on $G$ is lazy then 
\(\tau(\eps) = O \left( \log( |V(G)|/\eps) \,/\, \gap(P) \right).\)
\end{theorem}
\begin{theorem}\label{mixgap:theorem}(~\cite[Proposition 1.ii]{Si92} 
or \cite[Chapter 4]{AlFi})
For any time reversible Markov chain $P$ and $\eps>0$,
\(\gap(P) = \Omega(\log(1/2\eps)\,/\,\tau(\eps)).\)
\end{theorem}

\section{Composing simple permutations}

Another building block that we use are results on reversible computation
that enables us to compose simple permutations to construct permutations
that are easier to work with.  A classical result of Coppersmith and
Grossman~\cite{CoGr75} is that for $n > 3$ the set of width $2$ simple
permutations generates exactly the alternating group $A_n$.  Thus, all
compositions must be even permutations.

Formally, we define the set of width $w$ simple permutations$, \S_w$, 
as the set of permutations
$f_{i,J,h}$ where $i \in [n]$, $J = \{j_1,\ldots,j_w\}$ is a size $w$ ordered
subset of $[n]\setminus \{i\}$, and $h$ is a Boolean function on $\F2^w$. 
The permutation $f_{i,J,h}$ maps 
$(x_1,\ldots,x_n) \in \F2^n$ to 
$(x_1,\ldots,x_{i-1},x_i \oplus h(x_{j_1},\ldots,x_{j_w}),x_{i+1},\ldots,x_n)$.
We are primarily interested in width $2$ simple permutations, 
and denote $\S=\S_2$.

\begin{theorem} (Barenco et al.~\cite{BBCDMSSSW95})\label{wideand:theorem}
The permutation that flips the $n$-th bit of input $x$ if and only if the first
$w$ bits of $x$ are $1$ can be implemented as a composition
of $O(w)$ permutations from $\S$, as long as $w \leq n-2$.
\end{theorem}

\begin{theorem} (Brodsky~\cite{Br04})\label{basic-3-cycle:theorem}
for any distinct $x,y,z \in \F2^n$, one can compose $O(n)$
permutations from $\S$ to obtain the $3$-cycle $(xyz)$.
\end{theorem}

A length $\ell$ {\em implementation} of the permutation $\sigma$ is a sequence 
of permutations $\sigma_1,\ldots,\sigma_\ell$ from $\S$ whose composition 
is $\sigma$.
Theorem~\ref{basic-3-cycle:theorem} gives a length $O(n)$ implementation
for $3$-cycles.
We would like to use this implementation to construct a multicommodity flow
with low load on all edges.  However, Theorem~\ref{basic-3-cycle:theorem}
does not guarantee this. We solve this problem by randomizing the 
implementation, enabling us to prove a stronger theorem.

A length $\ell$ {\em randomized implementation} of the permutation
$\sigma$ is a sequence of {\em random} permutations
$\sigma_1,\ldots,\sigma_\ell$ from $\S$ whose composition is $\sigma$.
In Theorem~\ref{3-cycle:theorem} we give a randomized implementation
for 3-cycles, such that applying any prefix
$\sigma_1\cdots\sigma_{\ell'}$ of the randomized implementation of a
uniformly random 3-cycle $(xyz)$ to $x$ yields a string that looks
random. Namely, its {\em min-entropy} $H_\infty(\cdot)$ is high, which
is the minimum amount of information revealed when exposing the value
of a random variable $X$, that is
$H_\infty(X)=\min_\chi(-\log_2(\Pr[X=\chi]))$.

\begin{theorem}\label{3-cycle:theorem}
Let $x,y,z \in \F2^n$ be uniformly distributed and distinct. Then there
is a length $L=O(n)$ randomized implementation $\sigma_1 \cdots \sigma_L$
of the 3-cycle $(xyz)$ such that for all $\ell \in [L]$ 
the min-entropy of 
$(x\sigma_1\cdots\sigma_{\ell-1},\sigma_\ell)$ (which is a random variable on
$\F2^n \times\S$) is at least $\log_2( 2^n\cdot n^3 ) - O(1)$.
\end{theorem}
Note, this implies that the min-entropy of the marginals is big, i.e., 
$H_\infty(x\sigma_1\cdots\sigma_{\ell-1}) \geq n - O(1)$
and $H_\infty(\sigma_\ell) \geq \log_2(n^3) - O(1)$.

\section{Proof of Theorem~\ref{k2n3:theorem}}\label{k2n3proof:section}

In order to prove that the composition of random permutations from $\S$ 
approaches $k$-wise independence quickly we construct the Schreier graph 
$G_{k,n}=\schreier(\S,X^{(k)})$, where $X^{(k)}$ is the set of
$k$-tuple with $k$ distinct elements from the base set $X = \F2^n$.
It is convenient to think of $X^{(k)}$ as the set of $k$ by $n$ binary
matrices with distinct rows. 
A simple permutation acts on $X^{(k)}$ by acting on each of the rows. 
Then $P$ is the transition matrix of the random walk on $G_{k,n}$.
We prove that the random walk on this graph mixes rapidly.

To prove that $1/\gap(P) = O( n^2 k )$,
we first observe that $\gap(P)$ is monotone nonincreasing in $k$. 
This follows from the fact that the graph $G_{k+1,n}$ is a lift
of $G_{k,n}$ and therefore inherits the spectrum of $G_{k,n}$.
To see this, observe that any eigenfunction of $G_{k,n}$, can be lifted
to an eigenfunction on $G_{k+1,n}$, where the value of the latter on 
some $k+1$ by $n$ matrix is the value of the former on the matrix obtained by
deleting the last row. The eigenvalues of these two eigenfunctions is the 
same. 
In light of this observation, it is sufficient to prove the following two 
lemmas:

\begin{lemma}\label{cayley-gap:lemma}
$1/\gap(P) = O(n^2\cdot 2^n)$ for $k=2^n-2$.
\end{lemma}

\begin{lemma}\label{n2k-gap:lemma}
$1/\gap(P) = O(n^2 k)$ for $k \leq 2^n/3$.
\end{lemma}

We obtain the lower bound on the spectral gap of $P$ using the comparison 
technique~\cite{DiSa93}. This technique enables one to lower bound
$\gap(P)$ by $\gap(\tP)/A$, where $\tP$ is some other Markov chain,
and $A$ is the comparison constant. 
In our case, all chains are walks on regular graphs.
An upper bound on $A$ is obtained by constructing
a multicommodity flow on the underlying graph of $P$.
The flow flows a unit between
all pairs of endpoints of edges of $\tP$
such that the flow through each edge of $P$ is small.
To prove Lemmas~\ref{cayley-gap:lemma} and \ref{n2k-gap:lemma},
we compare $P$ to two different Markov chains. 
We start with the first Lemma.

\begin{proof}(of Lemma~\ref{cayley-gap:lemma})
For $k = 2^n - 2$, the state space of $P$ comprises
all even permutations of $\F2^n$. 
Let $\tP$ be a Markov chain on this state space, where
in each step 
we pick three distinct elements of the cube $x,y,z \in \F2^n$
and perform the permutation $(x y z)$.
It follows from a result of Friedman~\cite{Fr00}, 
that $1/\gap(\tP) = \Theta(2^n)$.
Therefore, it is sufficient to prove that the comparison constant of $P$ to
$\tP$ is $O(n^2)$.
\footnote{
Alternately, one can define a transition of $\tP$ as performing two random 
transpositions (not necessarily disjoint) and use a result of Diaconis and 
Shahshahani~\cite{DiSh81} that $1/\gap(\tP) = \Theta(2^n)$.}

To bound the comparison constant $A$, we need to construct a multicommodity 
flow $f$ in $G_{k,n}$ that flows a unit between every two matrices $M,M'$ 
such that $\tP(M,M')>0$.
Since the chains $P$ and $\tP$ correspond to random walks on regular graphs 
with degrees $d=\Theta(n^3)$ and 
$\tilde{d} = \Theta(2^{3n})$ respectively, 
the formula given in~\cite[Theorem 2.3]{DiSa93} reduces to:
\begin{eqnarray}\label{Abound:eqn}
A = 
(d/\tilde{d}) \cdot
\max_{(N,N') \in E(G_{k,n})} \left\{
\sum_{\gamma :\: (N,N') \in \gamma} f(\gamma) \cdot |\gamma|
\right\}.
\end{eqnarray}

Let $M,M'$ be two matrices such that $\tP(M,M')>0$.  Then $M'$ can be
obtained by applying some $3$-cycle $(xyz)$ to $M$.  Recall that the
randomized implementation given by Theorem~\ref{3-cycle:theorem}
induces a probability distribution on the length $L$ sequences of
permutations from $\Sigma$ whose composition is $(xyz)$.  Such a
distribution naturally translates to a distribution on length $L$
paths from $M$ to $M'$. We obtain a unit flow from $M$ to $M'$ by
flowing through each such path $\gamma$ an amount equal to its
probability.  We claim that the multicommodity flow obtained by
repeating this process for all $M,M'$ pairs satisfying $\tP(M,M')>0$
yields a small comparison constant.

Since $|\gamma| \cdot (d/\tilde{d}) = \Theta(n \cdot |\Sigma|/2^{3n})$
for all paths $\gamma$ with non-zero flow, the problem of bounding the
sum in (\ref{Abound:eqn}) reduces to bounding the total flow
through a given edge $e \in E(G_{k,n})$.  Let
$\gamma=(M_0,\ldots,M_L)$ be a path from $M_0$ to $M_L$, 
where $M_L$ is obtained from $M_0$ by applying the 3-cycle $(xyz)$. 
Assume further that $\gamma$ goes
through the edge $e$ at the $\ell$-th step, and that $x$ is
the $r$-th row of $M$. For any of the $\Theta(2^{4n}\cdot n)$ possible
assignments to $x,y,z,\ell,r$, we can determine the distribution
of the $r$-th row of the matrices $M_0,\ldots,M_L$. In particular, the
probability that $(M_{\ell-1},M_L)$ is equal to $e$ is bounded by the
probability that they coincide in their $r$-th row. By
Theorem~\ref{3-cycle:theorem}, in average over all assignments to
$x,y,z,\ell,r$, this probability is $O(1/2^n|\Sigma|)$. 
Putting it all together yields that, up to a constant factor, the
comparison constant $A$ is bounded $(n\cdot |\Sigma|/2^{3n}) \cdot
(2^{4n}\cdot n) \cdot (1/2^n |\Sigma|) = n^2$, as claimed.
\end{proof}

\begin{proof}(of Lemma~\ref{n2k-gap:lemma})
Let $\tP$ be the a Markov chain on the same state space as $P$,
which is the $k$ by $n$ binary matrices with distinct rows.
If the current state of $\tP$ is the matrix $M$, 
then the next state is determined by picking
a row $r \in \{1,\ldots,k\}$ and setting it to a random new value that
is distinct from all other $k-1$ rows. 
The process $\tP$ is the Markov chain of coloring the clique on 
$k$ vertices with $2^n$ colors defined in~\cite[section 4.1]{Je03}. 
Proposition 4.5 therein bounds its mixing time by 
$\tilde{\tau}(\eps) = O(k\log(k/\eps))$ as long as $k \leq 2^n/3$.
Setting $\eps=1/4k$ in Theorem~\ref{mixgap:theorem} implies that 
$\gap(\tP)=\Omega(1/k)$.
Therefore, as in the proof of Lemma~\ref{cayley-gap:lemma},
it is sufficient to prove that the comparison constant of $P$ to $\tP$ 
is $O(n^2)$.

Given matrices $M,M'$ such that $\tP(M,M')>0$, we know that $M'$ is obtained
from $M$ by changing the value of some row $r$ from $x$ to $y$.
To construct paths from $M$ to $M'$, we note that $M'$ can be obtained by 
applying the 3-cycle $(xyz)$ to $M$
for any $z \in \F2^n$ that is distinct from all rows of $M,M'$.
We choose $z$ at random from the $2^n-(k+1)$ allowed values.
As in the proof of Lemma~\ref{cayley-gap:lemma}, the randomized implementation
of $(xyz)$, given by Theorem~\ref{3-cycle:theorem}, defines 
a distribution on paths from $M$ to $M'$ and therefore a multicommodity
flow. We turn to bound the comparison constant, given by (\ref{Abound:eqn}).

As before, $|\gamma| \cdot (d/\tilde{d}) = \Theta(n \cdot |\Sigma|/k2^n)$
for all $\gamma$ with non-zero flow, and it suffices to bound the flow through
some edge $e \in E(G_{k,n})$. We enumerate over the choices of the 
position $\ell$, row $r$ and distinct $x,y$, which make a total of 
$\Theta(nk2^{2n})$ possible values. 
Again we apply Theorem~\ref{3-cycle:theorem} to argue that in average,
the probability of agreement with $e$ is bounded by $O(1/|\Sigma|2^n)$.
\footnote{
One should note that $z$ is uniformly distributed only over
$2^n-(k+1) > 2^{n-1}$ values. However, this is equivalent to conditioning a
uniform $z$ on an event with probability at least half and therefore 
(by Lemma \ref{cond:lemma}) can only increase the probability of agreement 
with $e$ by a factor of two.}
Therefore, up to a constant factor, 
$A=(n \cdot |\Sigma|/k2^n) \cdot (nk2^{2n}) \cdot (1/|\Sigma|2^n)=n^2$, 
as claimed.
\end{proof}

\section{Proof of Theorem~\ref{k2n2:theorem}}

In light of inequality (\ref{submultmix}), it is sufficient to
prove that $\tau(1/4) = \tO(n^2 k^2)$. 
The outline of the proof is the following.
We start by introducing the notion of a {\em generic} matrix,
and as suggested by the name, most matrices are generic.
The proof then proceeds by arguing that after a short random walk 
almost surely all matrices encountered are generic.
Therefore, it is sufficient to bound the mixing time of a walk that is 
restricted to generic matrices.
For such a walk, we can compare the chain to a chain defined {\em only} on 
generic matrices and achieve a much smaller comparison constant.
This yields the desired bound, $\tO(n^2 k^2)$.

Let $w=10\cdot(\log k + \log n)$. By assumption, we have $w \leq n/4$
for a sufficiently large $n$, and we set $p=\lceil n/2w \rceil$.
Let $C_1,\ldots,C_p,C$ be a partition of $[n]$ such that 
$|C_i|=w$ for $i=1,\ldots,p$ and $|C|=n-pw$. 
Consequently, $n/4 \leq n/2-w < |C| \leq n/2$.

We say that a $k$ by $n$ matrix is {\em generic}, if for all $j \in [p]$,
its restriction to $C_j$ has distinct rows.
It is not difficult to check that a uniformly distributed matrix $M$ is
almost surely generic. 
In fact, it is sufficient that the rows of $M$ are $2^{-w}$-close 
to $2$-wise independent, since then the probability that
$M$ is not generic is bounded by $p$ times the probability that the restriction
of $M$ to $C_j$ doesn't have distinct rows. 
This yields the bound
$p \cdot \binom{k}{2} \cdot (2 \cdot 2^{-w}) = o(1/n^3k^3)$ and
implies the following lemma:

\begin{lemma}\label{generic-2wise:lemma}
If the rows of a random $k$ by $n$ matrix $M$ are $2^{-w}$-close to $2$-wise 
independent, then $M$ is generic with probability $1-o(1/n^3k^3)$.
\end{lemma}

It follows from a result of Chung and Graham about the mixing time 
of the ``Aldous Cube''~\cite{ChGr97}, that the number of steps needed
to come close to $2$-wise independence, which is the same as the mixing time of 
$G_{2,n}$, is $O(n\log n)$.
This is stated in the following lemma (whose proof is deferred to 
Section~\ref{more:section}).

\begin{lemma}\label{two-wise:lemma}
For all $w\ge 1$ the $\eps$ mixing time of the Schreier graph
$\schreier(\S_w,X^{(2)})$ is $O(n \log n\log(1/\eps))$.
\end{lemma}

Therefore, the matrix obtained after 
$T_1=O(n \log n \cdot w) = O(n\log n \cdot(\log k + \log n))$
steps is $2^{-w}$-close to $2$-wise independent,
and by Lemma~\ref{generic-2wise:lemma} it is
generic with probability $1-o(1/n^3k^3)$.
This implies that if we proceed by $T_2 = O(n^3k^3)$ steps, 
then all $T_2$ matrices encountered are generic
with probability $1-o(T_2/n^3k^3) > 1-\eps_1$, for any fixed $\eps_1>0$ and 
sufficiently large $n$.

We introduce a new Markov chain $P'$. 
The state space of $P'$ consists of all generic $k$ by $n$ matrices. 
If the chain is currently at the matrix $M$, then the next state is determined
as follows. We pick a uniformly distributed simple permutation $\sigma \in \S$.
If $M\sigma$ is generic, we move to $M\sigma$. Otherwise, we remain at $M$.
Let $\tau'(\eps)$ denote the $\eps$-mixing time of $P'$,
and require that $T_2 \geq \tau'(\eps_2)$ for some fixed $\eps_2>0$.

We claim that as long as $2\eps_1+\eps_2 < 1/4$ the mixing time of $P$ can be 
bounded by $\tau(1/4) \leq T_1+T_2$.
To see this, 
let $M$ be some $k$ by $n$ matrix with distinct rows, and consider following 
two matrices.
The first matrix $M'$ obtained when starting at $M$ and walking $T_1+T_2$ 
steps using $P$.
The second matrix $M''$ is defined as follows. 
Let $\hat{M}$ be the matrix obtained
when starting at $M$ and performing $T_1$ steps of $P$.
If $\hat{M}$ is not generic, we set $M''=\hat{M}$.
Otherwise, $M''$ is the matrix reached by the length $T_2$ walk using $P'$
that starts at $\hat{M}$.
We claim that $d(M',M'') \leq \eps_1$ and that $M''$ is 
$(\eps_1+\eps_2)$-close to the uniform distribution over $k$ by $n$ matrices 
with distinct rows
\footnote{Note that by our assumptions, the distance between the uniform 
distribution over matrices with distinct rows and generic matrices is $o(1)$}.
Proving those claims will imply that 
\begin{eqnarray}
\tau(1/4) \leq \tau'(\eps_2) + O(n\log n \cdot(\log k + \log n)),
\end{eqnarray}
as long as $\tau'(\eps_2) = O(n^3k^3)$.

We start by checking that indeed $d(M',M'') \leq \eps_1$. It is convenient to 
think of the two length $T_1+T_2$ walks from $M$ to $M'$ and $M''$ as
defined over the same probability space $\S^{T_1+T_2}$ which is the choice
of a simple permutation in each of the $T_1+T_2$ steps. Denote the
the $P$-walk by $(M_0=M,M_1,\ldots,M_{T_1+T_2}=M')$.
Then, if all the matrices $M_{T_1},\ldots,M_{T_1+T_2}$ are generic,
it coincides with the walk leading to $M''$, and in particular we have
$M'=M''$. By the previous arguments, this event happens at least
with probability $1-\eps_1$, implying that $d(M',M'') \leq \eps_1$.

The proof that $M''$ is $(\eps_1+\eps_2)$-close to uniform is more delicate.
We know that the matrix $\hat{M}$ is generic with probability at least 
$1-\eps_1$.
Also, since $T_2 \geq \tau'(\eps_2)$, we know that conditioned on $\hat{M}$
being generic, $M''$ is $\eps_2$-close to the uniform distribution.
Therefore $M''$ is $(\eps_1+\eps_2)$-close to the uniform distribution over 
matrices with distinct rows. This argument can be easily formalized using 
Lemma~\ref{twoeps:lemma} of Section~\ref{more:section}.

We are left with the proof of the following lemma.

\begin{lemma}\label{genericmixing:lemma}
$\tau'(1/4) = \tO(n^2k^2)$.
\end{lemma}

To bound the mixing time of the Markov chain $P'$, we apply the comparison 
technique~\cite{DiSa93}. We compare $P'$ to the Markov chain $\tP$ defined
on the same state space, the $k$ by $n$ generic matrices. 
Given that $\tP$ is at a matrix $M$, we determine the next state as follows.
With probability half we pick a random column $c \in C$ and row $r \in [k]$
and flip the corresponding bit with probability half.
Otherwise, we pick at random an index $i \in [p]$, a row $r \in [k]$ and 
a string $\alpha \in \F2^w$ that is distinct from all other $k-1$ rows in 
the restriction of $M$ to the columns $C_i$. We set the bits at row $r$ 
and columns $C_i$ to $\alpha$.

Consequently, the following two lemmas,
imply Lemma~\ref{genericmixing:lemma}. Note that we need not worry about the
smallest eigenvalue of $P'$ since a random permutation from $\S$ is the 
identity with probability $1/16$.

\begin{lemma}\label{spectral:lemma}
$\gap(\tP) = \Omega(1/nk)$.
\end{lemma}

\begin{lemma}\label{simulation:lemma}
The comparison constant $A$ of $\tP$ to $P'$ satisfies $A = \tO(1)$.
\end{lemma}

\begin{proof} (of Lemma~\ref{spectral:lemma})

Consider two Markov chains $\tP_1$ and $\tP_2$:
\begin{enumerate}
\item
The state space of $\tP_1$ are the $k$ by $w$ binary matrices with distinct 
rows. At each step one chooses a random row and sets it to a random new
value distinct from all other $k-1$ rows. 
This chain is exactly the coloring chain of a clique on $k$ vertices
with $2^w$ colors of~\cite[Proposition 4.5]{Je03},
and as in the proof of Lemma~\ref{n2k-gap:lemma}, 
it satisfies $\gap(\tP_1)=\Omega(1/k)$.
\item
$\tP_2$ is the random walk on the $(n-wp)\cdot k$ dimensional binary cube, 
where in each step with probability half, 
one flips a random coordinate. Therefore,  $\gap(\tP_2)=\Omega(1/nk)$.
\end{enumerate}
One can think of the chain $\tP$ as the product of $p$ copies of $\tP_1$ and
one copy of $\tP_2$. Indeed the state space of $\tP$ is the direct product of
the $p+1$ state spaces. Moreover, a step of $\tP$ performs a move of $\tP_2$ 
with probability $1/2$ and otherwise performs the move in a randomly selected 
copy of $\tP_1$. 
It is straight forward to check that the spectral gap of $\tP$ is
$\min(\gap(\tP_1)/p,\gap(\tP_2))/2$, implying the desired bound.

\end{proof}

\begin{proof}(of Lemma~\ref{simulation:lemma})

Let $G'$ be the underlying graph of $P'$.
The vertices of $G'$ are the generic $k$ by $n$ matrices, 
and $(N,N')$ is an edge of $G'$ if $P'(N,N')>0$.
To bound the comparison constant $A$, we need to construct a multicommodity 
flow $f$ in $G'$ that flows a unit between every two matrices $M,M'$ 
such that $\tilde{P}(M,M')>0$.
The chains $P'$ and $\tP$ correspond to random walks on regular graphs
with degrees $d'=\Theta(n^3)$, $\tilde{d} = \Theta(kn2^w/w)$ respectively,
and as before the comparison constant $A$ is defined by (\ref{Abound:eqn}).

To build a path $\gamma$ from $M$ to $M'$ we need to distinguish two types
of $\tP$ transitions.
Type (i) flips the bit at row $r$ and column $c \in C$.
Type (ii) changes the bits at row $r$ and columns 
$C_i$ from $\alpha$ to $\alpha'$.
We start by constructing the type (i) paths.

Let $j \in [p]$ be a random index, and let $\beta \in \F2^w$ be the 
restriction of the $r$-th row of $M$ to $C_j$.
Also let $S$ be a random sequence of $w-1$ distinct elements from 
$C \setminus\{c\}$.
The unit flow from $M$ to $M'$ is along paths $\gamma=\gamma_{M,M'}^{S,j}$.
Each such path is defined by composing simple permutations from $\S$
to achieve the permutation that acts on $x \in \F2^n$ by flipping coordinate 
$c$ if the restriction of $x$ to $C_j$ is $\beta$. 
Clearly such a permutation maps $M$ to $M'$.
We follow the method of Barenco et al.~\cite{BBCDMSSSW95} to build 
an AND gate with $w$ inputs.
This gate inverts its output bit (the coordinate $c$) if its $w$ inputs (the
coordinates $C_j$) have some fixed value $\beta$. 
The coordinates in the set $S$ are used as ``scratch''.

Let $C_j=\{j_1,\ldots,j_w\}$, $S=\{s_1,\ldots,s_{w-1}\}$ and 
$\beta=(b_1,\ldots,b_w)$.
Let $\sigma_1$ be the simple permutation that flips coordinate 
$s_1$ of $x \in \F2^n$ if $x_{j_1}$ is equal to $b_1$,
and let $\sigma_\ell$ for $2 \leq \ell \leq w-1$ be the simple permutation 
that flips coordinate $s_\ell$ if $x_{s_{\ell-1}}$ is one and $x_{j_\ell}$ 
is equal to $b_\ell$. Also, we denote by $\tau_c$ the simple permutation that 
flips $x_c$ if $x_{s_{w-1}}$ is one and $x_{j_w}$ is equal to $b_w$. 
We claim that the following permutation
flips coordinate $c$ of $x \in \F2^n$ if the restriction of $x$ to $C_j$ is 
equal to $\beta$: 
\begin{eqnarray*}
\sigma = 
(\tau_c \sigma_{w-1} \cdots \sigma_2\sigma_1\sigma_2 \cdots \sigma_{w-1}
)^2
\end{eqnarray*}
To see this, one checks by induction that 
$\sigma_\ell\cdots\sigma_1\cdots\sigma_\ell$ flips coordinate $s_\ell$ if 
$x_{j_1}, \ldots, x_{j_\ell}$ is equal to $b_1,\ldots,b_\ell$.

For the type (ii) paths, we need to change the bits at row $r$ and
columns $C_i$ from $\alpha$ to $\alpha'$. 
The problem is that if we change $\alpha$ to $\alpha'$ bit by bit,
as suggested by the construction of type (i) paths, 
we might violate row distinctness.
To solve this problem, we start our path by applying a length 
$L = O(w\log w\cdot(1+2\log k))$ sequence $\phi$ of simple permutations with 
indices restricted to $C_i$.
Let $\hM = M\phi$ and $\hpM = M'\phi$, and let $C_i'$ and $C_i''$ be the first
and last $\lfloor (w-1)/2 \rfloor$ columns of $C_i$.
We say that $\phi$ is valid if for both the restriction of $\hM$ to $C_i''$ and
for the restriction of $\hpM$ to $C_i'$, have distinct rows. 
By Lemma~\ref{two-wise:lemma} we know for a random $\phi$, both $\hM$
and $\hpM$ are $1/8k^2$-close to $2$-wise independence. 
Therefore, a random $\phi$ is not valid with probability bounded by 
$k^2\cdot(2^{-w/2+1}+1/8k^2) \leq 1/4$.
If $\phi$ is valid we define a path $\gamma=\gamma_{M,M'}^{S,j,\phi}$
from $M$ to $M'$, where $j \in [p]\setminus \{i\}$ and $S$ is 
a length $w-1$ sequence of elements from $C$. 
The path is prefixed by $\phi$ to get from $M$ to $\hM$
and is suffixed by $\phi^{-1}$ to get from $\hpM$ to $M'$.
Let $\halpha$ and $\hpalpha$ be the restriction of the $r$-th row 
of $\hM$ and $\hpM$ to $C_i$ respectively,
and let $\beta$ be the restriction of the $r$-th row of $M$ to $C_j$.
Then the middle path connecting $\hM$ to $\hpM$ is defined as follows:
\begin{eqnarray*}
\sigma = [(
\prod_{\{c \in C_i'\,:\,\halpha_c \neq \hpalpha_c\}} \tau_c
) \cdot \sigma_{w-1} \cdots \sigma_2\sigma_1\sigma_2 \cdots \sigma_{w-1} ]^2
\cdot
[(
\prod_{\{c \in C_i \setminus C_i'\,:\,\halpha_c \neq \hpalpha_c\}} \tau_c
) \cdot \sigma_{w-1} \cdots \sigma_2\sigma_1\sigma_2 \cdots \sigma_{w-1} ]^2,
\end{eqnarray*}
where $\tau_c$ and $\sigma_\ell$ are as defined for the type (i) sequences.
Therefore it is guaranteed that the matrices encountered 
along the first and second half of the sequence agree with 
$\hM$ on the columns $C_i''$ and with 
$\hpM$ on the columns $C_i'$ respectively. 
Since $\phi$ is valid, this implies that we never attempt to move to 
a non-generic matrix throughout the entire path.
We define the unit flow from $M$ to $M'$ by splitting the flow uniformly
between all valid paths $\gamma$ designated by $S,j,\phi$. 

There are two points that need special attention in the constructed type (i) 
and type (ii) paths. 
The first point is that all indices of the simple permutations used in $\phi$
are in $C_i$. This is unacceptable for us, as it induces an undue load on
a small subset of $\S$. To solve this problem we replace
each simple permutation used in $\phi$ by a constant length sequence that 
avoids that problem. 
For example, the permutation that flips coordinate $i_1$ if $i_2$ and $i_3$ 
are $1$, denoted $\chi_{i_1,i_2,i_3}$, is replaced by the sequence
\( (\chi_{s_2,s_1,i_3} \chi_{s_1,i_2}, \chi_{s_2,s_1,i_3}, \chi_{i_1,s_2})^2\)
where permutation $\chi_{i_1,i_2}$ XORs coordinate $i_1$ with $i_2$.

The second point is that some of the simple permutations used 
($\sigma_1$ and some of the permutations in $\phi$) 
do not use three indices. However, in the definition of $\S$, 
we have three indices at our disposal even if we don't use all three.
We use this to guarantee that all simple permutations used have one index in 
$C_j$ and two from $S$ or $c$ for type (i) paths or $C_i$ for type (ii) paths.

To complete the proof, we have to bound the comparison constant $A$
given by (\ref{Abound:eqn}). 
We have $d'/\tilde{d} = \theta(n^2w/k2^w)$ and $|\gamma|=O(L)$.
Also, $f(\gamma)$ is $\Theta(w/n(m)_{w-1})$ for type (i) paths 
and $\Theta(w/|\S^{(w)}|^Ln(m)_{w-1})$ for type (ii) paths, 
where we denote $m=|C|$, 
$(m)_q = m(m-1)(m-2)\cdots(m-q+1)$,
and $\S^{(w)}$ as the width $2$ simple permutations restricted to the 
$w$-dimensional cube.
Therefore, we only have to bound the maximal number of $\gamma_{M,M'}^{S,j}$
and $\gamma_{M,M'}^{S,j,\phi}$ paths through an edge $(N,N')$.

We start with type (i) paths.
The first step is to extract as much information as possible about a path 
$\gamma$ through $(N,N')$ by considering the simple permutation $s$ associated 
with $(N,N')$.
Note first that $s$ determines $j$. 
Moreover, since only one of $\sigma_1,\ldots,\sigma_{w-1}$ and 
$\tau_c$ can be equal to $s$, any path $\gamma$ using $s$, 
must use it in one of $O(1)$ possible positions. 
Since a permutation $\sigma_\ell$ for $\ell \in [w-1]$ or $\tau_c$ determines 
two indices of $S,c$ 
there are only $\Theta((m)_{w-2})$ choices for $S,c$ that are 
consistent with $s$.
The last thing still needed to reconstruct $\gamma$ is the string 
$\beta \in \F2^w$. Since the columns $C_j$ are not modified throughout the
entire sequence, $\beta$ must be the restriction of some row of $N$ to $C_j$,
limiting $\beta$ to one of $k$ possible values.
Therefore, the total number of type (i) paths through $(N,N')$ 
is $O(k \cdot (m)_{w-2})$, 
and the contribution of the type (i) sequences to $A$ is:
\begin{eqnarray*}
A_{(i)} 
= O( \overset{d'/\tilde{d}}{\overbrace{(n^2w/k2^w)}} 
\cdot \overset{f(\gamma)\cdot|\gamma|}{\overbrace{(Lw/(m)_w)}}
\cdot \overset{\mbox{choices for }j,S,c,\beta\mbox{ and position}}
              {\overbrace{(k \cdot (m)_{w-2})}} )
= O(Lw^2/2^w) = o(1).
\end{eqnarray*}

For type (ii) paths we distinguish the cases where $(N,N')$ is 
in the first middle or last sections of a path 
$\gamma_{M,M'}^{S,j,\phi}$.
Consider the first section (and similarly the last).
We enumerate over possible positions $\ell \in [L]$. 
Then we know two indices of the
sequence $S$ and one of the $3L$ indices in $C_i$ that where used by $\phi$.
Therefore, we have $L\cdot(m)_{w-3}\cdot|\S^{(w)}|^L/w$ possible values for 
$S,i,\phi$ and the position.
This enables us to determine $M$ and $\hM$. We still have to determine the
row $r$, the two strings $\alpha, \alpha'$ and the index $j$
which have $O(kn2^w/w)$ possibilities.
Therefore the contribution of the first and last sections of type (ii)
paths is:
\begin{eqnarray*}
A_{(ii.\mbox{first,last})} 
&=& O( \overset{d'/\tilde{d}}{\overbrace{(n^2w/k2^w)}} 
\cdot \overset{f(\gamma)\cdot|\gamma|}{\overbrace{(Lw/|\S^{(w)}|^L(m)_w)}}
\cdot \overset
       {\mbox{choices for }j,S,i,\phi,\alpha,\alpha',\beta\mbox{ and position}}
       {\overbrace{(k2^wL\cdot(m)_{w-2}\cdot|\S^{(w)}|^L/w^2)}})\\
&=& O(L^2) = O(w^2\log^2 w\cdot(1+\log k)^2).
\end{eqnarray*}

For the middle section of type (ii) paths, as for the type (i) argument,
given $(N,N')$ we first determine the position up to $O(1)$ possible choices.
Then we determine the index $i$ or $j$ and two indices from $S$,
then we have $O( (m)_{w-2} \cdot |\S^{(w)}|^L/w )$ possibilities for 
$i,j,S,\phi$.
Also we have $k2^w$ choices for the row and the strings $\beta,\alpha$ and 
$\alpha'$.
Therefore,
\begin{eqnarray*}
A_{(ii.\mbox{middle})}
= O(\overset{d'/\tilde{d}}{\overbrace{(n^2w/k2^w)}} 
\cdot \overset{f(\gamma)\cdot|\gamma|}{\overbrace{(Lw/|\S^{(w)}|^L(m)_w)}}
\cdot \overset
       {\mbox{choices for }j,S,i,\phi,\alpha,\alpha',\beta\mbox{ and position}}
       {\overbrace{k2^w \cdot (m)_{w-2} \cdot |\S^{(w)}|^L/w}}
)
= O(Lw).
\end{eqnarray*}

\end{proof}


\section{Proof of Theorem~\ref{3-cycle:theorem}}

\newcommand{\Core}{{\rho_\mathit{core}}}
\newcommand{\Top}{\mathit{top}}
\newcommand{\Bot}{\mathit{bot}}
\newcommand{\cI}{\I}
\newcommand{\cIp}{{\I}^\prime}
\newcommand{\vp}{v^\prime}
\newcommand{\vpp}{v^{\prime\prime}}
\newcommand{\up}{u^\prime}
\newcommand{\tvp}{\tilde{v}^\prime}
\newcommand{\tvpp}{\tilde{v}^{\prime\prime}}

First, we describe the randomized implementation of a $3$-cycle
$(xyz)$ using the simple permutations in $\Sigma$.  Second, we show
that this randomized implementation satisfies the statement of the
theorem. The randomness is introduced into the implementation of
$(xyz)$ by using a permutation $\phi \in S_n$ and two vectors $v_4,v_5$.

Let $\phi$ be some permutation of the $n$ coordinates. 
If $\omega=\sigma_1\cdots\sigma_L$ implements $(x\phi,y\phi,z\phi)$, 
then $\omega^\phi$ is an implementation of $(xyz)$, 
where $\omega^\phi=\phi \omega \phi^{-1}$ is the conjugation of
$\omega$ with $\phi$, i.e. the conjugation each of the permutations
$\sigma_i$ used in $\omega$. Note that the set $\Sigma$ of simple
permutations is closed under conjugation by permutations from $S_n$, because
this just relabels the indices.

For a vector $v \in \F2^n$, we denote the first $n-2$ bits of $v$ by $\vp
\in \F2^{n-2}$ and the last two bits of $v$ by $\vpp \in \F2^2$, i.e.,
$v = \vp\vpp$.  We call the last two bits the {\it control bits}.  
For convenience, the notation
$\vp00$, $\vp01$, $\vp10$, and $\vp11$ denotes bit vectors comprising
the first $n-2$ bits of $v$ and the control bits $00$, $01$, $10$, and
$11$, respectively.  Let $(v)_j$ denote the $j$-th bit of a vector $v$.
Finally, let $v_1 = x\phi$, $v_2 = y\phi$ and $v_3 = z\phi$.

If $\vp_1$ is equal to $\vp_2$ or to $\vp_3$ then we say that $\phi$ is
invalid. 
This can only occur if $x$, $y$, or $z$ are less than Hamming distance
$3$ apart and $\phi$ maps all indices on which $x$
and $y$ (or $z$) differ to the control indices.  
For the rest of the description we assume that $\phi$ is valid.
Let $v_4, v_5 \in \F2^n$ be two additional vectors satisfying the validity
requirement of being at least Hamming distance $3$ from each other and from 
the former three vectors. 

Observe that $(v_1,v_2,v_3) = \psi_1 \psi_2$ where 
$\psi_1=(v_1,v_2)(v_4,v_5)$ and $\psi_2=(v_1,v_3)(v_4,v_5)$.
Therefore it suffices to implement the two double transpositions
$\psi_1$ and $\psi_2$. These are implemented in an identical manner.
Each implementation is divided into 15 blocks: a core block, which
implements the permutation $\Core = (\vp_500,\vp_501) (\vp_510,\vp_511)$,
and seven block pairs conjugating it.

The first four of these blocks, called $\pi$-blocks ensure that the
control bits of each of the four vectors are distinct.  Specifically, $\vp_i
\vpp_i$ is mapped to $\vp_i c_i$, where $c_1 = 00$, $c_2 = c_3 = 01$, $c_4
= 10$ and $c_5 = 11$.  If $\vpp_i = c_i$ then the corresponding block,
labeled $\pi_i$ performs a nop.  Otherwise, block $\pi_i$ performs the
permutation $(\vp_i \vpp_i, \vp_i c_i) (\vp_i a_i, \vp_i b_i)$ where
$\{a_i,b_i\} = \F2^2 \backslash \{\vpp_i,c_i\}$.

The remaining three blocks, called $\tau$-blocks, map $\vp_1$, $\vp_2$ (or
$\vp_3$), and $\vp_4$ to $\vp_5$, using the control bits to distinguish
between the four vectors.  Block $\tau_i$ performs the permutation
$\tau_i = \prod_{\vp \in \F2^{n-2}} (\vp c_i, \up c_i)$, where $\up =
\vp \oplus \vp_i \oplus \vp_5$.  Since it can easily be checked that
$\tau_i = \tau_i^{-1}$, that $\pi_i = \pi_i^{-1}$, and that
\[\pi_1 \pi_2 \pi_4 \pi_5 \tau_1 \tau_2 \tau_4 \Core 
  \tau_4 \tau_2 \tau_1 \pi_5 \pi_4 \pi_2 \pi_1 = \psi_1
\mathrm{\ \ \ and\ \ \ }
  \pi_1 \pi_3 \pi_4 \pi_5 \tau_1 \tau_3 \tau_4 \Core 
  \tau_4 \tau_3 \tau_1 \pi_5 \pi_4 \pi_3 \pi_1 = \psi_2,\]
we need only describe the implementation of each of these blocks.

Each of the blocks is implemented using $O(n)$ simple permutations.  Each
$\tau$-block is implemented by concatenating $n-2$ simple permutations,
where for $j = 1 \cdots n-2$, the $j$-th simple permutation is the
identity if $(\vp_i)_j = (\vp_5)_j$, and otherwise flips the $j$-th bit of
vector $v$ if $\vpp = c_i$.

The implementation of the $\Core$ and $\pi$ blocks is more involved.
Permutation $\Core$ flips bit $(\vpp)_2$ if and only if $\vp = \vp_5$.
Barenco et al~\cite{BBCDMSSSW95} showed how such permutations can be
implemented using $O(n)$ simple permutations, comprising four sub-blocks:
$\rho_\Top\rho_\Bot\rho_\Top\rho_\Bot$ where permutation $\rho_\Top$ flips
bit $(\vpp)_1$ if the first $\lceil (n-2)/2 \rceil$ bits of $\vp$ match
the first $\lceil (n-2)/2 \rceil$ bits of $\vp_5$, and where permutation
$\rho_\Bot$ flips bit $(\vpp)_2$ if the latter $\lfloor (n-2)/2 \rfloor$
bits of $\vp$ match the latter $\lfloor (n-2)/2 \rfloor$ bits of $\vp_5$
and $(\vpp)_1 = 1$.  Each sub-block uses the remaining $\lceil (n-2)/2
\rceil$ bits as ``scratch'', returning them to their original state by
the end of the sub-block.  For details about the construction of the
two sub-blocks see~\cite{BBCDMSSSW95} or Lemma~\ref{simulation:lemma}.

Each block $\pi_i$ is implemented in a similar manner using two
permutations that are nearly identical to the implementation of $\Core$.
The first (second) permutation performs the identity if $(\vpp)_1 =
(c_i)_1$ (respectively, $(\vpp)_2 = (c_i)_2$) and otherwise flips bit
$(\vpp)_1$ (respectively, $(\vpp)_2$) if $\vp = \vp_i$.

The length of the implementations of $\psi_1$ and $\psi_2$ is $O(n)$,
since each of the seven blocks can be implemented using $O(n)$ simple
permutations from $\Sigma$. The randomize implementation of $(xyz)$ is
obtained by uniformly choosing at random a valid permutation $\phi$ and
the two valid random vectors $v_4,v_5$.

We now prove that this randomized implementation satisfies the statement of the
theorem. 
Let $\Omega = \{x, y, z, v_4, v_5, \phi\}$ 
be the probability space obtained by uniformly
choosing three distinct vectors $x$, $y$, and $z$, and then uniformly
choosing a corresponding implementation, which is fixed by $v_4$, $v_5$,
and $\phi$.  Each point $\omega=(x,y,z,v_4,v_5,\phi)\in\Omega$
corresponds to an implementation $\sigma_1\cdots\sigma_L$ of the 
3-cycle $(xyz)$.
The size of $\Omega$ is $\Theta(2^{5n}n!)$, and   
although not uniform, the probability
of each point in $\Omega$ is $O(1/2^{5n}n!)$.  Thus, our problem
of upper-bounding 
$\Pr[x\sigma_1\sigma_2 \cdots \sigma_{\ell-1} = \tilde{x},\,
     \sigma_\ell = \tilde{\sigma}]$
reduces to a counting problem.

For all implementations $\omega\in\Omega$, the indices of the $\ell$-th 
permutation $\sigma_\ell$ depend {\em only} on its position, $\ell$, and $\phi$. 
Moreover, as we change $\phi$ the indices of the $\ell$-th permutation of the 
implementation of $ (x,y,z,v_4,v_5,\phi)$ agree with the indices of some 
fixed permutation $\tilde{\sigma}$ only on a subset of $S_n$ that is of 
size $O(n!/n^3)$ and depends only on $\ell$ and $\tilde{\sigma}$.

To establish the theorem we need to prove that for any given $\phi$ the 
number of choices of $x$, $y$, $z$, $v_4$, and $v_5$, such that 
$x\sigma_1\sigma_2\cdots\sigma_{\ell-1} = \tilde{x}$, is $O(2^{4n})$, implying
the number of points in $\Omega$ that agree with $\tilde{x}$
and $\tilde{\sigma}$ is $O(2^{4n}n!/n^3) = O(|\Omega|/2^n n^3)$.  This is 
accomplished by the following lemma:

\begin{lemma}\label{counting:lemma}
Let $\phi \in S_n$ be fixed.
Then the set of all $x,y,z,v_4,v_5$ such that implementation
corresponding to $(x,y,z,v_4,v_5,\phi)$ satisfies the equality
$x \sigma_1\sigma_2 \cdots \sigma_{\ell-1} = \tilde{x}$
is of size $O(2^{4n})$.
\end{lemma}

\begin{proof}(of Lemma~\ref{counting:lemma})

Let $v_1 = x\phi$, $v_2 = y\phi$, $v_3 = z\phi$, and $\tilde{v}=\tilde{x}\phi$.
Let $\Omega_{\tilde{v},\ell}$ be the set of tuples $(v_1,\ldots,v_5)$ for 
which $x \sigma_1\sigma_2 \cdots \sigma_{\ell-1} = \tilde{x}$ is satisfied. 
Note that this set is independent of $\phi$.
Then the claim is that $|\Omega_{\tilde{v},\ell}|=O(2^{4n})$.

The proof is via case analysis with respect to position $\ell$.
Without loss of generality we assume that the position is in the first
half of the implementation, that which realizes permutation $(v_1,v_2)
(v_4,v_5)$, otherwise, swapping $v_2$ and $v_3$ allows the same argument
to be reused for the latter half of the implementation.  Furthermore,
due to symmetry, we assume that the position of $\ell$ is in or to the
left of block $\Core$.  There are four main cases: either $\ell$ is on
a boundary between two blocks, $\ell$ is in block $\tau_i$, $\ell$ is in the
block $\Core$, or $\ell$ is in block $\pi_i$.

\begin{figure}[ht]
\large
\[\underbrace{v_1 \stackrel{\pi_1}{\longrightarrow} 
        \vp_1 c_1 \stackrel{\pi_2\pi_4\pi_5}{-\!\!\!\!-\!\!\!\!\longrightarrow}
        \vp_1 c_1}_{\tvp \mathrm{\ fixes\ }\vp_1}
                  \stackrel{\tau_1}{\longrightarrow}
  \underbrace{\vp_5 c_1 \stackrel{\tau_2\tau_4}{\longrightarrow} 
        \vp_5 c_1 \stackrel{\Core}{\longrightarrow} 
        \vp_5 c_2 }_{\tvp \mathrm{\ fixes\ }\vp_5}
        \cdots
\]
\caption{The evolution of $v_1$. \label{fig:thm9}}
\end{figure}

In the first case, the position, $\ell$, is on a block boundary.
Since each $\pi$-block only toggles bits $(\vpp)_1$ and $(\vpp)_2$,
if position $\ell$ is adjacent to a $\pi$-block, then $\tvp = \vp_1$.
Thus, all but two bits of $v_1$ are fixed by $\tilde{v}$.  If position,
$\ell$, is on a boundary but is not adjacent to a $\pi$-block, then it
must occur after block $\tau_1$.  Since block $\tau_1$ maps $\vp_100$
to $\vp_500$, and none of the remaining blocks, $\tau_i$ or $\Core$,
change the $\vp$ component to any other value, we have $\tvp = \vp_5$.
Thus, all but two bits of $v_5$ are fixed by $\tilde{v}$, implying that
$|\Omega_{\tilde{v},\ell}| = O(2^{4n})$.

In the second case, the position, $\ell$, is inside block $\tau_i$.  If $i
\not= 1$, then none of the simple permutations in block $\tau_i$ flips a
bit.  Therefore, the value of $\tvp = \vp_5$; thus fixing all but two bits
of $v_5$, as before.  If $i = 1$, then at position $\ell$, we know exactly
how many of the $n-2$ simple permutations have already been performed.
Let $j$ be this number.  Hence we know that $\tvp = (\vp_5)_1, \ldots,
(\vp_5)_j, (\vp_1)_{j+1}, \ldots, (\vp_1)_{n-2}$.  Therefore, $j$ bits of
$v_5$ and $n - 2 - j$ bits of $v_1$ are therefore fixed by $\tilde{v}$,
implying that $|\Omega_{\tilde{v},\ell}| = O(2^{4n})$ as well.

In the third case, the position, $\ell$, is inside block $\Core$. In this
case we must look at the sub-blocks of the block $\Core$.  If the position
occurs on a sub-block boundary, and since each of the sub-blocks simply
toggles the bits $(\vpp)_1$ and $(\vpp)_2$, the remaining bits of $\vp_5$
are fixed by $\tvp$.  If the position $\ell$ is inside a sub-block, then
things are only slightly more complicated.  Assume that position $\ell$
is in a $\rho_\Top$ sub-block (similar arguments hold for $\rho_\Bot$).
Then, $\rho_\Top$ toggles bit $(\vpp)_1$ if the first half of $\vp$
matches the first half of $\vp_5$.  The bits being matched are never
modified and the other half of the bits of $\vp$ are used as ``scratch''.
We know that the first half of $\vp$ and $\vp_5$ coincide throughout the
block $\rho_\Top$, and therefore $\tvp$ determines this half of $\vp_5$.
The operations on the ``scratch'' half depends only on the fixed half and
the position, and therefore can be reversed,  reducing the problem to the
position occurring at the beginning of $\rho_\Top$.  Thus, $\tilde{v}$
fixes all but two of the bits of $v_5$.

In the last case, the position, $\ell$, is inside a $\pi$-block.  Block
$\pi_i$ comprises two blocks that are similar to $\Core$.  Each of the
two blocks is either the identity or toggles $(\vpp)_1$ or $(\vpp)_2$ if
$\vp = \vp_i$.  If $i = 1$ then the two blocks in Block $\pi_i$ behave
in the same manner as block $\Core$, except that $\tilde{v}$ fixes all
but two of the bits of $v_1$ rather than $v_5$.  If $i \not = 1$, then,
for the most part, the argument remains the same.  We need only consider
what happens if the position, $\ell$, is in one of the eight sub-blocks.
As mentioned before, half of the bits of $\vp$ are not modified by the
sub-block, while the other half are used as ``scratch''.  Again, without
loss of generality, we assume that the sub-block does not modify the
first half of $\vp$.  As before, $\tvp$ fixes the first half of $\vp_1$.
We enumerate on all choices for the first half of $\vp_i$.  This enables
us to reverse the operations of the sub-block on the ``scratch'', fixing
the second half of $\vp_1$---as in the third case.  This implies that
$|\Omega_{\tilde{v},\ell}| = O(2^{4n})$, and completes the proof.

\end{proof}

\section{Odds and Ends}\label{more:section}

\begin{proof}(of Lemma~\ref{two-wise:lemma})\\
We have to prove that for all $w\ge 1$ the mixing time of
$G_{2,n}^{(w)}=\schreier(\S_w,X^{(2)})$ is $O(n \log n)$.

Given a $2$ by $n$ matrix with rows $s,t$, 
we change basis to $s,u$ with $u=s \oplus t$. 
Let $i\in[n]$ be a random coordinate, and consider the action of a width 
$w$ permutations XORing the $i$-th bit with a random function $h$ on $w$ 
distinct coordinates from $[n]\setminus\{i\}$.
We claim that its action on $s,u$ is the same as XORing the $i$-th bit of 
$s$ and $u$ with two {\em independent} random bits 
$\alpha_s$ and $\alpha_u$ respectively.
The bits $\alpha_s,\alpha_u$ are one with probability
$1/2$ and $p_\ell = 1 - \prod_{j=1}^w (1-\frac{\ell}{n-j})$ respectively,
where $\ell$ is the number of ones in $u$ not counting the $i$-th bit.
To see that this is indeed the resulting walk we observe the fact that
if $s$ and $t$ differ on one of the input bits of the random function $h$, 
then the value of the $i$-th coordinate of $s$ and of $t$ change
independently with probability half. Otherwise they change simultaneously
with probability $1/2$.

The $u$-component of this walk is a variant of the Aldous cube, and by
the comment at the end of~\cite{ChGr97} it follows that this walk
mixes in $O(n \log n)$ time. We are left to show that in this time the
walk on both components mixes. The way to see it is to notice that in
$O(n \log n)$ time the event $A$ where the indices $i$ assume all
possible values in $1,2,\ldots,n$ (coupon collector) happens with high
probability. Now since the bits $\alpha_s$ are independent of
$\alpha_u$, we get that even when we condition over the walk on the
$u$ component, the $s$ component achieves uniform distribution
conditioned on $A$, which ends the proof.
\end{proof}

\begin{lemma}\label{twoeps:lemma}
Let $A$ be an event such that $\Pr[A] \geq 1-\eps$,
and let $Z$ be a random variable over a domain $\Omega$ such
that $d(Z|A,\mbox{uniform}) \leq \eps$.
Then $d(Z,\mbox{uniform}) \leq 2\eps$.
\end{lemma}
\begin{proof}
\[d(Z,\mbox{uniform})
=
\max_{S \subseteq \Omega} \Pr[Z \in S]-\frac{|S|}{|\Omega|}
\leq
\max_{S \subseteq \Omega} \Pr[Z \in S|A]+\Pr[\overline{A}]-\frac{|S|}{|\Omega|}
\leq
\eps + d(Z|A,\mbox{uniform}) \leq 2\eps.
\]
\end{proof}

\begin{lemma}\label{cond:lemma}
Let $X$ be a random variable and $A$ an event. 
Then $\Pr[X|A] \leq \Pr[X]/\Pr[A]$.
(Follows from the definition of conditional probability.)
\end{lemma}

\section{Some concluding remarks}\label{conclude:section}

Let us review what we currently know about the spectral gap of the Markov chain
$P=P_\Sigma^{(k,n)}$. 
By Theorem~\ref{k2n3:theorem}, $\gap(P) \leq \Omega(1/n^2k)$. 
On the other hand, $\gap(P)$ is nonincreasing in $k$ by the
lifting argument from Section~\ref{k2n3proof:section}.
Since for $k=1$, $P$ is the standard random walk on the cube,
we have that $\gap(P) \geq 1/n$.

In general, a generating set $S$ for which the spectral gap is large becomes 
more difficult as $k$ increases, until the largest conceivable $k$, which is
$2^n-2$. In this case, this is the random walk on the Cayley graph of the 
alternating group $A_N$ for $N=2^n$ with the generating set $S$.
It is open whether one can find a constant size set 
for which $A_N$ is an expander, \cite[Problem 10.3.4]{Lu94}.
\footnote{The problem of finding a constant size expanding set for 
$A_N$ or $S_N$ is equivalent.}
On the other hand, by Alon and Roichman~\cite{AlRo94}, a random set 
of permutations of size $O(N \cdot \log N)$ will almost surely have a constant 
spectral gap. 
Although smaller expanding sets for $A_N$ are not known to exist,
the general belief is that such sets exist;
Rozenman, Shalev, and Wigderson assume the existence of an 
$N^{1/30}$ expanding set for $A_N$,~\cite[section 1.4]{RSW04}.

Our results suggest that width $2$ permutations may be used to
construct an $O(\log^3 N)$ expanding set for $A_N$. 
However, several obstacles stand in the way of achieving this goal.
The first one is to prove that for width $2$ permutations the spectral gap
does not deteriorate with $k$, as we believe, and is $\Omega(1/n)$ for all $k$.
The second problem is to achieve a constant gap.  
To this end, one has to overcome the inherent and obvious weakness of the 
width $2$ simple permutations.
Namely, that their action depends only on two 
coordinates and changes only one.  
This leads to poor expansion because there is only a small chance that the 
action will flip a specific bit
or increase the distance between two similar vectors.
One approach to avoiding this problem is to replace the standard set
of generators of the cube $e_1,\ldots,e_n$ with some expanding set
of size $O(n)$. Such an expanding set for the cube can readily be
constructed from the generating matrix of a good code~\cite{DeSo91},
and could then be used to define
an $O(n^3)$ expanding set of permutations. 


\bibliography{bib}

\end{document}